\documentclass[11pt,a4paper]{scrartcl}
\usepackage{amsmath,amsfonts,amstext, amsthm, pgf, tikz, paralist, framed, color, colortbl, array, graphicx, setspace, cite, scrpage2}

\author{Marie-Louise Bruner \\ Institut f\"ur Diskrete Mathematik und Geometrie \\ Technische Universit\"at Wien \\ Vienna, Austria \\ \texttt{marie.louise-bruner@tuwien.ac.at}}
\title{A short note on the Stanley-Wilf Conjecture \\ for permutations on multisets\thanks{This work was supported by the Austrian Science Foundation FWF, grant S9608-N13.}}
\date{\today}

\theoremstyle{plain}
\newtheorem{Thm} {Theorem} [section]

\newtheorem{Def} [Thm]{Definition}

\theoremstyle{definition}
\newtheorem{Ex} [Thm]{Example}
\theoremstyle{remark}
\newtheorem{Rem} [Thm]{Remark}

\newcommand{\pe}{permutation}
\newcommand{\mul}{multiset}

\begin{document}

\maketitle

\begin{abstract}
The concept of pattern avoidance respectively containment in \pe s can be extended to \pe s on \mul s in a straightforward way. In this note we present an alternative proof of the already known fact that the well-known Stanley-Wilf Conjecture, stating that the number of \pe s avoiding a given pattern does not grow faster than exponentially, also holds for \pe s on \mul s.
\end{abstract}

\section{Introduction}

Since 1985 when Rodica Simion and Frank Schmidt published the first systematic study of \textit{Restricted Permutations} \cite{simion1985restricted} the area of pattern avoidance in \pe s has become a rapidly growing field of enumerative combinatorics.
A comprehensive overview of this area can be found in \cite{bona_combinatorics_2004, kitaev2011patterns}. We first state the central definition for \pe s on a set:

\begin{Def}
A \pe\ $p=p_1 p_2 \ldots p_n$ of length $n \geq k$ is said to \textit{contain} another \pe\ $q=q_1 q_2 \ldots q_k$ as a \textit{pattern } if we can find $k$ entries $p_{i_1}, p_{i_2}, ..., p_{i_k}$ with $i_1 < i_2 < ... < i_k$ such that $p_{i_a} < p_{i_b} \Leftrightarrow q_a < q_b$. In other words, \textbf{$p$ contains $q$} if we can find a subsequence of $p$ that is order-isomorphic to $q$. If there is no such subsequence we say that $p$ \textbf{avoids the pattern $q$}.
\label{def_avoidance_pes}
\end{Def}

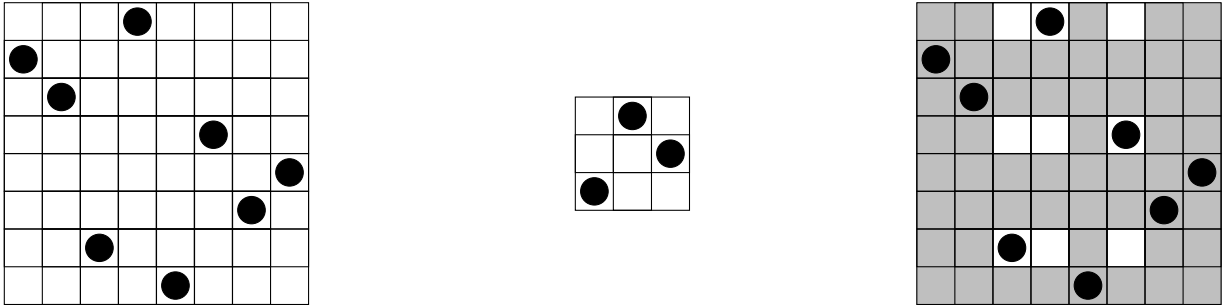
\begin{figure}
\begin{minipage}[hbt]{5cm}
	\centering
	\begin{tikzpicture}
	[0/.style={rectangle, draw, minimum size=28.5pt}, 
	1/.style={circle, draw, fill=black, minimum size=1.5pt}, 
	2/.style={rectangle, fill=black!25, minimum size=14.25pt}, 
	scale=0.5]

	\foreach \x/\name in {0, 1, 2, 3, 4, 5, 6}
    	\node[0] at (\x, 0) {};
	\foreach \x/\name in {0, 1, 2, 3, 4, 5, 6}
    	\node[0] at (\x, 1) {};
	\foreach \x/\name in {0, 1, 2, 3, 4, 5, 6}
    	\node[0] at (\x, 2) {};    
	\foreach \x/\name in {0, 1, 2, 3, 4, 5, 6}
    	\node[0] at (\x, 3) {};    
	\foreach \x/\name in {0, 1, 2, 3, 4, 5, 6}
    	\node[0] at (\x, 4) {};    
	\foreach \x/\name in {0, 1, 2, 3, 4, 5, 6}
    	\node[0] at (\x, 5) {};    
	\foreach \x/\name in {0, 1, 2, 3, 4, 5, 6}
    	\node[0] at (\x, 6) {};    
    
	\foreach \x/\y in { 0/6, 1/5, 2/1, 3/7, 4/0, 5/4, 6/2, 7/3}
    	\node[1] at (-0.5+\x, -0.5+\y) {};
	\end{tikzpicture}
\end{minipage}
\hfill
\begin{minipage}[hbt]{1cm}
	\centering
	\begin{tikzpicture}
	[0/.style={rectangle, draw, minimum size=28.5pt}, 
	1/.style={circle, draw, fill=black, minimum size=1.5pt}, 
	2/.style={rectangle, fill=black!25, minimum size=14.25pt}, 
	scale=0.5]

	\foreach \x/\name in {0, 1}
    	\node[0] at (\x, 0) {};
	\foreach \x/\name in {0, 1}
    	\node[0] at (\x, 1) {};
       
	\foreach \x/\y in { 0/0, 1/2, 2/1}
    	\node[1] at (-0.5+\x, -0.5+\y) {};
	\end{tikzpicture}
\end{minipage}
\hfill
\begin{minipage}[hbt]{5cm}
	\centering
	\begin{tikzpicture}
	[0/.style={rectangle, draw, minimum size=28.5pt}, 
	1/.style={circle, draw, fill=black, minimum size=1.5pt}, 
	2/.style={rectangle, fill=black!25, minimum size=14.25pt}, 
	scale=0.5]

	\foreach \x in {0, 1, 2, 3, 4, 5, 6, 7}
    	\node[2] at (-0.5+\x, -0.5) {};
	\foreach \x in {0, 1, 2, 3, 4, 5, 6, 7}
    	\node[2] at (-0.5+\x, 1.5) {};
	\foreach \x in {0, 1, 2, 3, 4, 5, 6, 7}
    	\node[2] at (-0.5+\x, 2.5) {};
	\foreach \x in {0, 1, 2, 3, 4, 5, 6, 7}
    	\node[2] at (-0.5+\x, 4.5) {};
	\foreach \x in {0, 1, 2, 3, 4, 5, 6, 7}
    	\node[2] at (-0.5+\x, 5.5) {};
	\foreach \y in {0, 1, 2, 3, 4, 5, 6, 7}
    	\node[2] at (-0.5, -0.5+\y) {};    
	\foreach \y in {0, 1, 2, 3, 4, 5, 6, 7}
    	\node[2] at (0.5, -0.5+\y) {};    
	\foreach \y in {0, 1, 2, 3, 4, 5, 6, 7}
    	\node[2] at (3.5, -0.5+\y) {};
	\foreach \y in {0, 1, 2, 3, 4, 5, 6, 7}
    	\node[2] at (5.5, -0.5+\y) {};    
	\foreach \y in {0, 1, 2, 3, 4, 5, 6, 7}
    	\node[2] at (6.5, -0.5+\y) {};    

	\foreach \x/\name in {0, 1, 2, 3, 4, 5, 6}
    	\node[0] at (\x, 0) {};
	\foreach \x/\name in {0, 1, 2, 3, 4, 5, 6}
    	\node[0] at (\x, 1) {};
	\foreach \x/\name in {0, 1, 2, 3, 4, 5, 6}
    	\node[0] at (\x, 2) {};    
	\foreach \x/\name in {0, 1, 2, 3, 4, 5, 6}
    	\node[0] at (\x, 3) {};    
	\foreach \x/\name in {0, 1, 2, 3, 4, 5, 6}
    	\node[0] at (\x, 4) {};    
	\foreach \x/\name in {0, 1, 2, 3, 4, 5, 6}
    	\node[0] at (\x, 5) {};    
	\foreach \x/\name in {0, 1, 2, 3, 4, 5, 6}
    	\node[0] at (\x, 6) {};    
    
	\foreach \x/\y in { 0/6, 1/5, 2/1, 3/7, 4/0, 5/4, 6/2, 7/3}
    	\node[1] at (-0.5+\x, -0.5+\y) {};
    \end{tikzpicture}
\end{minipage}
\caption{The \pe\ $23718465$ (left) contains the pattern $312$ (middle), as can be seen by deleting the rows and columns marked in grey (right).}
\label{pattern_avoidance_grid}
\end{figure}

\begin{Ex}
The \pe\ $p=23718465$ (written in one-line representation respectively as a sequence of integers) contains the pattern $312$, since the entries $714$ (or several other examples) form a $312$-pattern. See Figure \ref{pattern_avoidance_grid} for a graphical representation with the help of \pe\ matrices, i.e. square binary matrices with exactly one entry equal to $1$ - placed at $(p(i),i)$ - per row and per column. As a matter of fact, this \pe\ contains all possible patterns of length three. This is not the case for patterns of length four: $p$ contains the pattern $2134$ as is shown by the entries $3146$ but $p$ avoids the pattern $4321$ since it contains no decreasing subsequence of length four. 
\end{Ex}

Considering \textit{\pe s on \mul s} (or words on a certain alphabet) is a natural generalization of ordinary \pe s. By a \mul\ we mean a ``set'', where elements may occur more than once. We write $\left\lbrace 1^{m_1}, 2^{m_2}, \ldots, n^{m_n}\right\rbrace $ for the \mul\ in which the element $i$ occurs $m_i$-times and call $m_i$ the \textit{multiplicity} of $i$. A \mul\ is called \textit{regular} in case $m_i=m$ for some integer $m$ and we use the notation $[n]_m$ for $\left\lbrace 1^{m}, 2^{m}, \ldots, n^{m}\right\rbrace $ where $[n]:=[n]_1$.
The definition of pattern avoidance respectively containment introduced in Definition \ref{def_avoidance_pes} can be extended to \pe s on \mul s in a straightforward way. 
For a certain pattern to be contained in a \pe , repetitions in the pattern have to be represented by repetitions in the \pe . For instance, the \mul -\pe\ $1214324$ contains the patterns $122$, $123$ and $321$ but avoids the pattern $211$. enumerative questions for restricted \pe s on \mul s have been studied during the last decade, see e.g. the work of Heubach and Mansour \cite{heubach2006avoiding} and Albert et. al. \cite{albert2001permutations}.

When attempting to compute the number $S_n(q)$ of ordinary \pe s of length $n$ avoiding a certain pattern $q$, very satisfying results were obtained for patterns of length three: $S_n(q)=c_n= \frac{1}{n+1}\binom{2n}{n}$, the $n$-th Catalan number, for all six $3$-patterns $q$ (see e.g. \cite{knuth1968art,simion1985restricted}). For longer patterns however, the situation turns out to be a lot more complicated: even for patterns of length four there is still one of the three \textit{Wilf-equivalence} classes\footnote{Two patterns $q_1$ and $q_2$ are \textit{Wilf-equivalent} if $S_n(q_1)=S_n(q_2)$ holds for all integers $n$.} for which no enumeration formula has yet been found. There is nevertheless a more general result on restricted \pe s, stating that the number of $n$-\pe s avoiding an arbitrary given pattern does not grow faster than exponentially:

\begin{Thm}
[Stanley-Wilf Conjecture, 1990, proven 2004 in \cite{klazar2000füredi,marcus2004excluded}] Let $q$ be an arbitrary pattern. Then there exists a constant\footnote{As a matter of fact, the constant $c_q$ only depends on the length of the forbidden pattern $q$.} $c_q$ such that for all positive integers it holds that 
\begin{eqnarray}
S_n(q) \leq c_q^n.
\end{eqnarray}
\label{Conj_Stanley_Wilf}
\end{Thm}

As Bona remarked in \cite{bona_combinatorics_2004}, this was quite an ambitious conjecture since it postulates that the number of $q$-avoiding $n$-\pe s does not grow faster than exponentially whereas the total number of $n$-\pe s grows super-exponentially, cf. the Stirling formula that states that $n!$ is asymptotically equal to $\sqrt{2 \pi n}  \left(n/\mathrm e\right)^{n}$.

A generalization of the Stanley-Wilf conjecture has been considered independently by Klazar and Marcus \cite{klazar2007extensions}  and by Balogh, Bollob\'as and Morris \cite{DBLP:journals/ejc/BaloghBM06}. Klazar and Marcus proved an exponential bound on the number of hypergraphs avoiding a fixed permutation, settling various conjectures of Klazar as well as a conjecture of Br\"and\'en and Mansour. Balogh, Bollob\'as and Morris went even further in their generalization and showed similar results for the growth of hereditary properties of partitions, ordered graphs and ordered hypergraphs. For details, please consider the original work. The results in \cite{DBLP:journals/ejc/BaloghBM06,klazar2007extensions} being very general however required rather involved proofs. In both papers, a generalized version of the F\"uredi-Hajnal conjecture (Theorem \ref{Conj_Fueredi_Hajnal}) was first formulated and the proof of several intermediary results was necessary.

From the results in both these papers it follows that the Stanley-Wilf conjecture also holds for permutations on multisets since these can be represented with the help of bipartite graphs which are (very) special cases of hypergraphs. To the best of our knowledge it has however not yet been stated in the literature that the Stanley-Wilf conjecture also holds for permutations on multisets. We therefore wish to stress this point here, as pattern avoidance in permutations on multisets has attracted a great deal of interest in the past few years. 

This present note is devoted to providing a simple and direct proof of the Stanley-Wilf Conjecture for \pe s on \mul s, given the fact that the Stanley-Wilf conjecture holds for ordinary permutations. The following proof uses the idea of Klazar \cite{klazar2000füredi} and does not require the employment of any other results or a further generalization of already known results.

\section{A generalization of Stanley-Wilf to \mul s}
\subsection{Multiset-Stanley-Wilf}

Our goal is to prove the following generalization of Theorem \ref{Conj_Stanley_Wilf} to \pe s on regular \mul s:
\begin{Thm}
[Multiset-Stanley-Wilf] Let $q$ be an arbitrary \pe\ on an ordinary set and $S_{n,m}(q)$ the number of \pe s on the regular \mul\ $[n]_m$ avoiding this pattern. Then there exists a constant $e_q$ merely depending on the length of $q$ so that the following holds for all positive integers $n$ and $m$:
\begin{eqnarray}
S_{n,m}(q) \leq e_q^{n\cdot m}.
\label{eqn_multi_sw}
\end{eqnarray}
\label{Multi_Stanley_Wilf}
\end{Thm}

Multiset-Stanley-Wilf can easily be extended to \pe s on an arbitrary \mul , see Remark \ref{rem_arbitrary_multiset}.

\begin{Rem}
Now, why would the upper bound suggested in \eqref{eqn_multi_sw} provide an equally strong result for \pe s on \mul s as Stanley-Wilf for ordinary \pe s? How does $e_q^{n\cdot m}$ compare to the total number of \pe s on the regular \mul\ $[n]_m$? For a regular \mul\ the total number of \pe s is equal to
\[
A_{n,m}:=\dbinom {m\cdot n}{\underbrace{m,m, \ldots, m}_{n-\text{times}}} = \dfrac{(mn)!}{(m!)^n}.
\]
Let us set $f(n,m):=\sqrt{2 \pi m n} \left( n^m / \left( \sqrt{2 \pi m} \cdot \mathrm e\right) \right)^n$. Applying Stirling's formula then yields that $f(n,m)~\ll~A_{n,m}$ when $m$ is fixed and $n \to \infty$. Obviously $f(n,m)$ grows a lot faster than $e_q^{n\cdot m}$ for any constants $e_q$ and $m$. 
\end{Rem}

\begin{Rem}
The condition that the avoided pattern $q$ is a \pe\ on an ordinary set is crucial in Theorem \ref{Multi_Stanley_Wilf}. Indeed, the number $S_{n,m}(q)$ may grow faster than exponentially if the pattern $q$ is a \mul -\pe . For example, consider $m$-Stirling \pe s, i.e. \pe s on a regular \mul\ with multiplicity $m$ avoiding the pattern $212$. As stated in \cite{kuba2010enumeration}, the number of $m$-Stirling \pe s with $n$ distinct elements is equal to $n!m^n \binom {n-1+1/m}{n}$ and thus grows super-exponentially.
\end{Rem} 

For the proof of this \mul -version of Stanley-Wilf we follow the proof of the original result. Stanley-Wilf was not proven directly but via another conjecture formulated by F\"uredi and Hajnal.  

\subsection{The F\"uredi-Hajnal Conjecture}
\label{subsec_FH}

In \cite{furedi1992davenport} a conjecture concerning pattern avoidance in binary matrices was presented. 

\begin{Def}
Let $P$ and $Q$ be matrices with entries in $\left\lbrace 0,1 \right\rbrace $ and let $Q$ have the dimension $m\times n$. We say that the matrix $P$ \textbf{contains} $Q$ as a pattern, if there is a submatrix $\tilde Q$ of $P$, so that $\tilde{Q}_{i,j}=1$ whenever $Q_{i,j}=1$ for $i \leq m$ and $j \leq n$. If there is no such submatrix $\tilde Q$, we say that $P$ \textbf{avoids} $Q$.
\label{def_avoidance_matrices}
\end{Def}

This means that $P$ contains $Q$ as a pattern, if, by deleting some rows and some columns, one can obtain a matrix $\tilde{Q}$ with the same size as $Q$ that has a $1$-entry everywhere where $Q$ has a $1$-entry. Note that $\tilde{Q}$ must not necessarily have its $0$-entries in the same places as in $Q$: $\tilde{Q}$ may have more $1$-entries than $Q$ but not less. 

\begin{Thm}[F\"uredi-Hajnal Conjecture]
Let $Q$ be any \pe\ matrix. We define $f(n,Q)$ as the maximal number of $1$-entries that a $Q$-avoiding $(n \times n)$-matrix $P$ may have. Then there exists a constant $d_Q$ so that
\[
f(n,Q) \leq d_Q \cdot n.
\]
\label{Conj_Fueredi_Hajnal}
\end{Thm} 

In the year 2000 Martin Klazar proved that the F\"uredi-Hajnal conjecture implies the Stanley-Wilf conjecture \cite{klazar2000füredi}. Four years later the F\"uredi-Hajnal conjecture was proven by Adam Marcus and G\'abor Tardos in \cite{marcus2004excluded}, finally providing a simple and ``gorgeous'' \cite{zeilberger2004how} proof of the long-standing Stanley-Wilf conjecture. In their proof Marcus and Tardos showed that the constant $c_q$ in Theorem \ref{Conj_Stanley_Wilf} does not depend on the pattern $q$ itself but merely on its length.

\subsection{Klazar's proof}

For the proof of Theorem \ref{Multi_Stanley_Wilf}, it will merely be necessary to show that F\"uredi-Hajnal implies the \mul -version of Stanley-Wilf. Let us therefore briefly recall the argument used by Martin Klazar in his proof. 

In order to establish a connection between pattern avoidance in matrices and pattern avoidance in \pe s Klazar takes an elegant detour via pattern avoidance in simple bipartite graphs and defines the following notion of pattern containment:

\begin{Def}
Let $P([n],[n'])$ and $Q([k],[k'])$ be simple bipartite graphs, where $k \leq n$ and $k' \leq n'$. Then we say that $P$ \textbf{contains $Q$ as an ordered subgraph} if two order preserving injections $f:[k]\rightarrow[n]$ and $f':[k']\rightarrow[n']$ can be found so that if $vv'$ is an edge of $Q$, then $f(v)f'(v')$ is an edge of $P$.
\end{Def}

Note again that - as in Definition \ref{def_avoidance_matrices} of pattern avoidance in binary matrices - the following holds: if $f(v)f'(v')$ is an edge of $P$, $vv'$ does not necessarily have to be an edge of $Q$. 

Clearly, every \pe\ can be identified with a simple bipartite graph in a unique way. For a \pe\ $p$ on $[n]$ the associated graph $G_p$ is the bipartite graph with vertex set $([n],[n])$ and where $e=(i,j)$ is an edge iff $p_i=j$ in $p$ ($i$ is an element of the first set of $n$ elements, $j$ is an element of the second one). Then the following is a direct consequence: If the \pe\ $p$ contains a \pe\ $q$ as a pattern, then $G_p$ contains $G_q$ as an ordered subgraph. Reversly, if $p$ avoids $q$, $G_p$ will also avoid $G_q$. However, not every simple bipartite graph corresponds to a \pe\ (this is only the case if the vertex degree is equal to one for all vertices), thus $S_n(q) \leq  G_n(q)$, where $G_n(q)$ is the number of simple bipartite graphs on $([n],[n])$ avoiding the graph $G_q$ corresponding to a \pe\ $q$. In his proof, Klazar therefore aims at showing that for every \pe\ $q$ there is a constant $c_q$ so that $G_n(q) \leq c_q^n$.

Let $P$ be a simple bipartite graph on $([n],[n])$ that avoids $G_q$. Then the adjacency matrix $A(P)$ of $P$ avoids the adjacency matrix $A(G_q)$ and Theorem \ref{Conj_Fueredi_Hajnal} implies that $A(P)$ can have at most $d_q \cdot n=d_{A(G_q)} \cdot n$ entries equal to $1$, resepectively that $G_p$ can have at most $d_q \cdot n$ edges. By gradually contracting the graph $P$ - reducing its size to half in every step without loosing the $G_q$-avoiding property - Klazar shows that this leaves at most an exponential number of possibilities for the graph $P$: $G_n(q) \leq 15^{2d_qn}$. Thus $S_n(q) \leq c_q^n$ with $c_q = 15^{2d_q}$.

\subsection{Using Klazar's idea for \pe s on \mul s}
We shall use the same idea of contracting simple bipartite graphs in order to prove Theorem \ref{Multi_Stanley_Wilf}. It is straightforward to see how simple bipartite graphs can be used to represent \pe s on \mul s. Let $p=p_1 p_2 \ldots p_l$ be a \pe\ on the \mul\ $[n]_m$, i.e. a \pe\ of length $l=n \cdot m$. Then the associated bipartite graph $G_p$ is a graph with vertex set $([l],[n])$ and where $e=(i,j)$, $i \in [l]$ and $j \in [n]$ is an edge iff $p_i=j$. Note that this bipartite graph is balanced iff $m=1$, i.e. in the case of \pe s on ordinary sets. See the left-hand side of Figure \ref{counterex_multiSW} for the representation of the \mul -\pe s $1212$ and $111$.

\begin{proof}[Proof of Theorem \ref{Multi_Stanley_Wilf}]
Let $G_{n,m}(q)$ be the number of simple bipartite graphs on $([n \cdot m],[n])$ avoiding the graph $G_q$ corresponding to the ordinary \pe\ $q$. We shall show that $G_{n,m}(q) \leq e_q^{n\cdot m} $, implying that $S_{n,m}(q) \leq e_q^{n\cdot m} $. As in the case of ordinary \pe s and balanced simple bipartite graphs it holds that $S_{n,m}(q) \leq G_{n,m}(q)$ since not all simple bipartite graphs on $([n \cdot m],[m])$ correspond to \pe s on the \mul\ $[n]_m$.

Let $G$ be a simple bipartite graph on $([n \cdot m],[n])$ avoiding $G_q$. In a first step we shall contract the bipartite graph $G$ to a balanced bipartite graph $G'$ on the vertex set $([n],[n])$. For this purpose we merge $m$ consecutive vertices of the $[n \cdot m]$-vertex set to a single one in the following way. If $i \in [n]$ and $j \in [n]$ are two vertices in $G'$, then let $(i,j)$ be an edge if there is at least one edge between the set of vertices $V_{m,i}:=\left\lbrace (m-1)i+1, (m-1)i+2, \ldots, mi \right\rbrace $ and $j$ in $G$. For an example of such a contraction, see Figure \ref{counterex_multiSW}.

The balanced bipartite graph $G'$ inherits the $G_q$-avoiding property from $G$. Indeed, if it was possible to find a ``copy'' of $G_q$ in $G'$, we could find ``ancestors'' of all the edges involved in this copy in $G$. Then these ``ancestor-edges'' would certainly form a copy of $G_q$ in $G$.

Let us recapitulate: The resulting graph $G'$ is a simple bipartite graph on $([n],[n])$ avoiding $G_q$. If we consider the adjacency matrices $A(G')$ respectively $A(G_q)$, i.e. the binary matrices for which $a_{i,j}=1$ iff $(i,j)$ is an edge in the corresponding graph, it is clear that $A(G')$ avoids $A(G_q)$. From F\"uredi-Hajnal (Theorem \ref{Conj_Fueredi_Hajnal}) it follows that $A(G')$ can have at most $d_{A(G_q)} \cdot n =d_q \cdot n$ entries equal to $1$, implying that $G'$ may have at most $d_q \cdot n$ edges. As remarked earlier in \ref{subsec_FH}, the constant $d_q$ merely depends on the length of $q$.

How many different graphs $G$ can lead to the same contracted graph $G'$? For every edge $(i,j)$ in $G'$, there are exactly $2^m-1$ different sets of edges between $V_{m,i}$ and $j$ leading to this specific edge. Indeed, every non-empty subset of $\left\lbrace(k,j): k \in V_{m,i} \right\rbrace $ leads to an edge between $i$ and $j$ in $G'$. Thus there are in total at most
\[
(2^m-1)^{d_q \cdot n}
\]
possible graphs $G$ contracting to the balanced graph $G'$ and we obtain:
\[
G_{n,m}(q) \leq (2^m-1)^{d_q \cdot n} \cdot G_{n,1}(q). 
\]

\begin{figure}[t]
\begin{center}
	\begin{tikzpicture}[scale=0.8]
	\tikzstyle{vertex}=[circle,fill=black!25,minimum size=7pt,inner sep=0pt]

	\foreach \x/\y/\name in {0/0/1, 2.5/0/2, 8/0/3, 10.5/0/4, 0/-1/5, 2.5/-1/6, 8/-1/7, 10.5/-1/8, 0/-2/9, 0/-3/10, 0/-5/11, 2.5/-5/12, 	8/-5/13, 10.5/-5/14, 0/-6/15, 0/-7/a}
    	\node[vertex] (A-\name) at (\x, \y) {};
    
	\foreach \from/\to in {1/2, 3/4, 3/8, 5/6, 7/4, 7/8, 9/2, 10/6, 11/12, 13/14, 15/12, a/12}
    	{ \draw[thick] (A-\from) -- (A-\to);  }

	\node at (3.5,-0.5) (n1) {};
	\node at (7,-0.5) (n2) {};
	\node at (5.25,-0.75) {contraction};
	\node at (3.5,-5.5) (n3) {};
	\node at (7,-5.5) (n4) {};
	\node at (5.25,-5.75) {contraction};
	\draw[->,very thick,black!35] (n1) -- (n2);
	\draw[->,very thick,black!35] (n3) -- (n4);

	\node at (1.25,0.5) {$1212$};
	\node at (1.25,-4.5) {$111$};    
	\end{tikzpicture}
\end{center}
\caption{Why the proof of Theorem \ref{Multi_Stanley_Wilf} does not work for \mul -patterns: The top-left graph avoids the bottom-left one but the contracted graph on the top-right no longer avoids the one on the bottom right.}
\label{counterex_multiSW}
\end{figure}
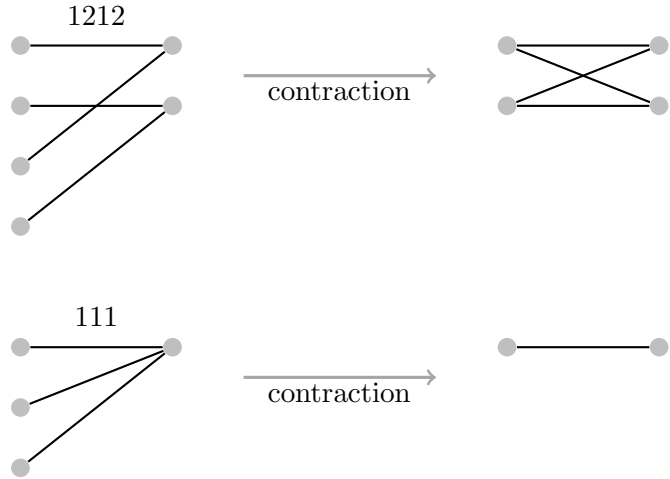

Now that we have obtained a balanced bipartite graph $G'$ avoiding $G_q$, we can continue in Klazar's proof where it was shown that
\[
G_{n,1}(q) = G_n(q) \leq 15^{2d_qn}.
\]
This finally leads to 
\[
G_{n,m}(q) \leq \left( (2^m-1)\cdot 15^2\right) ^{d_q \cdot n}, 
\] 
proving that 
$S_{n,m}(q) \leq e_q^{m \cdot n}$ where $e_q=\left(  2 \cdot 15^2 \right) ^{d_q}$ is a constant depending only on the length of the pattern $q$.
\end{proof}

\subsection{Concluding remarks}
\begin{Rem}
Using the proof given above, it is clear why the Stanley-Wilf conjecture cannot hold for permutations on multisets avoiding some multiset-pattern. Indeed, it can easily be seen why the contraction does not work if the forbidden pattern itself is a \pe\ on a \mul . The crucial point in our proof is that the balanced graph resulting from the contraction of the original bipartite graph $G$ still avoids the forbidden pattern-graph $G_q$. However, if the pattern itself must be contracted to a balanced bipartite graph, then the pattern-avoiding property is in general not passed on. For example, consider the (bipartite graph corresponding to the) \pe\ $p=1212$ represented in Figure \ref{counterex_multiSW}. It clearly avoids the (graph corresponding to the) pattern $q=111$, since no element occurs thrice in $p$. To the contrary, the contracted graph obviously contains the contracted pattern-graph since its edge-set is non-empty.
\label{rem_arbitrary_multiset}
\end{Rem}

\begin{Rem}
The above generalization of Stanley-Wilf to \mul s was formulated for regular \mul s - what happens if we consider arbitrary \mul s instead of regular ones? The reasoning and remarks above are obviously also valid if we set $m:=\max \left\lbrace m_i: i \in [n]\right\rbrace $ for an arbitrary \mul . The contraction done in the first step of the proof of Theorem \ref{Multi_Stanley_Wilf} must then  be carried out in the following way: merge the first $m_1$-many vertices, then the following $m_2$-many vertices and so on until finally the last $m_n$-many vertices are merged to a single one.
\end{Rem}

\section*{Acknowledgements}
I wish to thank Alois Panholzer for introducing me to the fascinating topic of restricted \pe s and for supporting my research.

\bibliographystyle{plain}
\bibliography{references}

\end{document}